\documentclass[12pt]{smfart}

\usepackage{smfenum}
\usepackage{smfthm}
\usepackage{amssymb}
\usepackage{mathrsfs}
\usepackage[francais]{babel}

\newcommand{\val}{\mathrm{val}}
\newcommand{\Qp}{\mathbf{Q}_p}
\newcommand{\Zp}{\mathbf{Z}_p}
\newcommand{\Fp}{\mathbf{F}_p}
\newcommand{\ZZ}{\mathbf{Z}}
\newcommand{\OO}{\mathcal{O}}
\newcommand{\MM}{\mathfrak{m}}
\newcommand{\Fpbar}{\overline{\mathbf{F}}_p}
\newcommand{\Qpbar}{\overline{\mathbf{Q}}_p}
\newcommand{\dcris}{\mathrm{D}_{\mathrm{cris}}}
\newcommand{\eps}{\varepsilon}
\renewcommand{\phi}{\varphi}
\renewcommand{\projlim}{\varprojlim}
\renewcommand{\injlim}{\varinjlim}
\newcommand{\Fil}{\mathrm{Fil}}
\renewcommand{\geq}{\geqslant}
\renewcommand{\leq}{\leqslant} 
\newcommand{\G}{\mathrm{GL}_2(\Qp)}
\newcommand{\B}{\mathrm{B}(\Qp)}
\newcommand{\K}{\mathrm{GL}_2(\Zp)}
\newcommand{\gal}{\mathcal{G}_{\Qp}}
\newcommand{\calA}{\mathscr{A}}
\newcommand{\calO}{\mathscr{O}_{\mathscr{E}}}
\newcommand{\pr}{\mathrm{pr}}
\newcommand{\ddiese}{\mathrm{D}^{\sharp}}
\newcommand{\dplus}{\mathrm{D}^+}
\newcommand{\dfont}{\mathrm{D}}
\newcommand{\Hom}{\mathrm{Hom}}
\newcommand{\LC}{\mathrm{LC}}

\newcommand{\ind}{\mathrm{ind}}
\newcommand{\Id}{\mathrm{Id}}
\newcommand{\res}{\mathrm{res}}
\newcommand{\Indbg}{\mathrm{Ind}_{\mathrm{B}}^{\mathrm{G}}}
\newcommand{\simB}{\underset{\mathrm{B}}{\sim}}
\newcommand{\simG}{\underset{\mathrm{G}}{\sim}}

\author[L. Berger]{Laurent Berger}
\address{CNRS \& IHES \\
Le Bois-Marie\\
35 route de Chartres\\
91440 Bures-sur-Yvette \\ 
France}
\email{laurent.berger@ihes.fr}
\urladdr{www.ihes.fr/\~{}lberger/}

\title[Repr\'esentations galoisiennes et repr\'esentations de $\mathrm{GL}_2(\mathbf{Q}_p)$]
{Repr\'esentations modulaires de $\mathrm{GL}_2(\mathbf{Q}_p)$ et repr\'esentations galoisiennes de dimension $2$} 

\alttitle{Modular representations of $\mathrm{GL}_2(\mathbf{Q}_p)$ and $2$-dimensional galois representations}

\date{Octobre 2005}


\NumberTheoremsIn{subsection}

\begin{document}

\begin{abstract}
On montre la conjecture de Breuil concernant la r\'eduction modulo $p$ des repr\'esentations triangulines $V$ et des repr\'esentations $\Pi(V)$ de $\mathrm{GL}_2(\mathbf{Q}_p)$ qui leur sont associ\'ees par la correspondance de Langlands $p$-adique. L'ingr\'edient principal de la d\'emonstration est l'\'etude de certaines repr\'esentations lisses irr\'eductibles de $\mathrm{B}(\mathbf{Q}_p)$ via des mod\`eles construits en utilisant les $(\varphi,\Gamma)$-modules.
\end{abstract}

\begin{altabstract}
We prove Breuil's conjecture concerning the reduction modulo $p$ of trianguline representations $V$ and of the representations $\Pi(V)$ of $\mathrm{GL}_2(\mathbf{Q}_p)$ associated to them by the $p$-adic Langlands correspondence. The main ingredient of the proof is the study of some smooth irreducible representations of $\mathrm{B}(\mathbf{Q}_p)$ through models built using the theory of $(\varphi,\Gamma)$-modules.
\end{altabstract}

\maketitle

\setcounter{tocdepth}{2}

\tableofcontents

\setlength{\baselineskip}{20pt}

\section*{Introduction}

Cet article s'inscrit dans le cadre de la correspondance de Langlands $p$-adique. Dans ses articles \cite{Br2,BRst}, Breuil a d\'efini une correspondance
qui \`a une repr\'esentation $p$-adique $V$ de dimension $2$ qui est cristalline ou semi-stable, associe une repr\'esentation $\Pi(V)$ de $\G$. Le fait que $\Pi(V)$ est non-nul, irr\'eductible et admissible a \'et\'e d\'emontr\'e par Colmez pour les repr\'esentations semi-stables dans \cite{CL} puis par Breuil et l'auteur pour les repr\'esentations cristallines dans \cite{BB}. La correspondance a ensuite \'et\'e \'etendue aux repr\'esentations triangulines par Colmez dans \cite{CT}. Dans \cite{Br1,Br2}, Breuil a par ailleurs d\'efini une correspondance en caract\'eristique $p$ et conjectur\'e qu'elle \'etait compatible avec la premi\`ere. L'objet de cet article est de d\'emontrer cette conjecture pour les repr\'esentations triangulines\footnote{Dans les cas dits {\og exceptionnels \fg}, c'est-\`a-dire $V$ cristalline non $\phi$-semi-simple ou bien $V$ non-cristalline mais le devenant sur une extension ab\'elienne, la d\'emonstration repose sur des r\'esultats non-r\'edig\'es.}.

Afin d'\'enoncer ce th\'eor\`eme, nous devons introduire quelques notations qui nous serviront par la suite. Dans tout cet article, $p$ est un nombre premier et le corps $L$ (le {\og corps des coefficients \fg}) est une extension finie de $\Qp$ dont on note $\OO_L$ l'anneau des entiers, $\MM_L$ l'id\'eal maximal et $k_L$ le corps r\'esiduel. On note $\gal$ le groupe de Galois $\mathrm{Gal}(\Qpbar/\Qp)$, $\omega$ la r\'eduction modulo $p$ du carac\`ere cyclotomique $\eps$, et pour $y \in L$ ou $y \in k_L$ on note $\mu_y$ le caract\`ere non-ramifi\'e de $\gal$ qui envoie $\mathrm{Frob}_p^{-1}$ (c'est-\`a-dire l'inverse du frobenius arithm\'etique) sur $y$. On normalise le corps de classes pour que $p\in\Qp^\times$ s'envoie sur $\mathrm{Frob}_p^{-1}$ ce qui permet de voir $\mu_y$ comme le caract\`ere non-ramifi\'e de $\Qp^\times$ qui envoie $p$ sur $y$. On note enfin $\omega_2$ le caract\`ere fondamental de Serre de niveau $2$. 

On appelle $\B$ le sous-groupe de Borel sup\'erieur de $\G$. Pour all\'eger les notations, il nous arrivera d'\'ecrire $\mathrm{G}$ pour $\G$ et $\mathrm{B}$ pour $\B$ ainsi que d'identifier $\Qp^\times$ au centre de $\G$. Si $\Pi_1$ et $\Pi_2$ sont deux $k_L$-repr\'esentations lisses de longueur finie de $\B$ (ou de $\G$), on \'ecrira $\Pi_1 \simB \Pi_2$ (ou bien $\Pi_1 \simG \Pi_2$) pour signifier que les semi-simplifi\'ees de $\Pi_1$ et de $\Pi_2$ sont isomorphes.

Si $V$ est une repr\'esentation $p$-adique de $\gal$, on en choisit un $\OO_L$-r\'eseau $T$ stable par $\gal$ et $\overline{V} = (k_L \otimes_{\OO_L} T)^{\mathrm{ss}}$ ne d\'epend pas du choix de $T$. De m\^eme, si $\Pi$ est une repr\'esentation de $\G$ qui admet un r\'eseau $\Pi^0$, on note $\overline{\Pi}$ la semi-simplifi\'ee de la r\'eduction de $\Pi^0$ qui, si elle est de longueur finie, ne d\'epend pas du choix du r\'eseau dans une classe de commensurabilit\'e.

Si $r \in \{0,\hdots,p-1\}$ et si $\chi : \Qp^\times \to k_L^\times$ est un caract\`ere continu, que l'on identifie \`a un caract\`ere continu de $\gal$ via le corps de classes, alors on pose :
\[ \rho(r,\chi) = (\ind_{\mathrm{Gal}(\Qpbar/\mathbf{Q}_{p^2})}^{\mathrm{Gal}(\Qpbar/\Qp)}(\omega_2^{r+1})) \otimes \chi. \]

On pose par ailleurs : 
\[ \pi(r,\lambda,\chi) = \left( \frac{\ind_{\K\Qp^\times}^{\G}
    \mathrm{Sym}^r k_L^2}{T-\lambda} \right) \otimes (\chi \circ
\det), \] 
o\`u $\lambda \in \Fpbar$ et $T$ est un certain op\'erateur de Hecke.

\begin{enonce*}{Th\'eor\`eme A}
Si $V$ est une repr\'esentation trianguline irr\'educible, alors :
\begin{align*}
\overline{V} = \rho(r,\chi) & \Leftrightarrow \overline{\Pi}(V) = \pi(r,0,\chi) \\
\overline{V} = \begin{pmatrix} \mu_\lambda \omega^{r+1} & 0 \\ 0 & \mu_{\lambda^{-1}} \end{pmatrix} \otimes \chi & \Leftrightarrow \overline{\Pi}(V) \simG \pi(r,\lambda,\chi) \oplus \pi([p-3-r], \lambda^{-1}, \omega^{r+1} \chi)
\end{align*}
\end{enonce*}

Ici $[p-3-r]$ est l'unique entier appartenant \`a $\{ 0, \hdots, p-2 \}$ et qui est congruent \`a $r$ modulo $p-1$.

La d\'emonstration de ce th\'eor\`eme est fond\'ee sur le lien entre $\Pi(V)$ et le $(\phi,\Gamma)$-module $\dfont(V)$, d'abord montr\'e par Colmez pour les repr\'esentations semi-stables puis par Breuil et l'auteur pour les repr\'esentations cristallines et enfin par Colmez pour les repr\'esentations triangulines (les cas exceptionnels n'\'etant pas (encore) r\'edig\'es).

Le premier chapitre est consacr\'e \`a des rappels. Le deuxi\`eme chapitre est consacr\'e \`a l'\'etude de repr\'esentations de $\B$ construites \`a partir de certaines repr\'esentations galoisiennes. On trouve soit des restrictions d'induites paraboliques, soit des restrictions de supersinguli\`eres et le fait d'en avoir des mod\`eles explicites est important pour la suite. Nous d\'emontrons notamment les r\'esultats suivants d'int\'er\^et ind\'ependant.

\begin{enonce*}{Th\'eor\`eme B}
Si $\Pi$ est une $k_L$-repr\'esentation lisse irr\'eductible de $\G$ admettant un caract\`ere central, alors :
\begin{enumerate}
\item si $\Pi$ est un caract\`ere, ou la sp\'eciale, ou supersinguli\`ere, alors sa restriction \`a $\B$ est toujours irr\'eductible;
\item si $\Pi$ est une s\'erie principale, alors sa restriction \`a $\B$ est une extension du caract\`ere induisant $\chi_1 \otimes \chi_2$ par une repr\'esentation irr\'eductible $\Indbg(\chi_1 \otimes \chi_2)_0$.
\end{enumerate}
\end{enonce*}

\begin{enonce*}{Th\'eor\`eme C}
Si $\Pi_1$ et $\Pi_2$ sont deux $k_L$-repr\'esentations lisses semi-simples et de longueur finie de $\G$ (dont les composantes irr\'eductibles admettent un caract\`ere central) et dont les semi-simplifications des restrictions \`a $\B$ sont isomorphes, alors $\Pi_1$ et $\Pi_2$ sont d\'ej\`a isomorphes en tant que repr\'esentations de $\G$.
\end{enonce*}

Dans le troisi\`eme chapitre, on utilise les r\'esultats explicites du deuxi\`eme pour d\'emontrer le th\'eor\`eme A. On termine en faisant le point sur les r\'esultats que l'on en d\'eduit concernant la r\'eduction modulo $p$ des repr\'esentations cristallines et semi-stables.

\section{Rappels et compl\'ements}\label{rappels}

L'objet de ce chapitre est de rappeler certaines des constructions faites dans le cadre de la {\og correspondance de Langlands $p$-adique \fg}. Nous renvoyons \`a \cite{Br1,Br2,BRst,BB,CL,CT} pour plus de d\'etails.

\Subsection{Repr\'esentations galoisiennes et $(\phi,\Gamma)$-modules}

Nous reprenons les notations de l'introduction. Les caract\`eres $\chi : \gal \to k_L^\times$ sont tous de la forme $\omega^r \mu_y$ pour $r \in \{0,\hdots,p-2\}$ et $y \in k_L^\times$. Rappelons la classification des repr\'esentations absolument irr\'eductibles de $\gal$ de dimension $2$ sur $k_L$. Soit $\omega_2 : \mathrm{Gal}(\Qpbar/\mathbf{Q}_{p^2}) \to \mathbf{F}_{p^2}^\times$ le caract\`ere fondamental de Serre de niveau $2$. On suppose d\'esormais que $\mathbf{F}_{p^2} \subset k_L$. Si $r \in \{0,\hdots,p-1\}$ et si $\chi : \Qp^\times \to k_L^\times$ est un caract\`ere continu, que l'on identifie \`a un caract\`ere continu de $\gal$ via le corps de classes, alors on pose :
\[ \rho(r,\chi) = (\ind_{\mathrm{Gal}(\Qpbar/\mathbf{Q}_{p^2})}^{\mathrm{Gal}(\Qpbar/\Qp)}(\omega_2^{r+1})) \otimes \chi. \]
On obtient ainsi toutes les repr\'esentations absolument irr\'eductibles de $\gal$ de dimension $2$ sur $k_L$, et les entrelacements entre les $\rho(r,\chi)$ sont les suivants : 
\[ \rho(r,\chi) \simeq  \rho(r,\chi \mu_{-1}) \simeq  \rho(p-1-r,\chi \omega^r)  \simeq \rho(p-1-r,\chi \omega^r \mu_{-1}). \]

La classification des repr\'esentations $p$-adiques de $\gal$ est, comme on le sait, beaucoup plus compliqu\'ee mais la th\'eorie des $(\phi,\Gamma)$-modules permet de s'y retrouver un petit peu. Soit $\Gamma = \mathrm{Gal}(\Qp(\mu_{p^\infty})/\Qp)$ et $\calO$ l'anneau
$\calO = \{ \sum_{i \in \ZZ} a_i X^i$ o\`u $a_i \in \OO_L$
et $a_{-i} \to 0$ quand $i \to \infty \}$. On munit cet anneau d'un
frobenius $\OO_L$-lin\'eaire $\phi$ d\'efini par $\phi(X)=(1+X)^p-1$
et d'une action $\OO_L$-lin\'eaire de $\Gamma$ donn\'ee par $\gamma(X)
= (1+X)^{\eps(\gamma)}-1$ si $\gamma \in \Gamma$. Un
$(\phi,\Gamma)$-module \'etale est un $\calO$-module $\dfont$ de type fini
muni d'un frobenius semi-lin\'eaire $\phi$ tel que $\phi^*(\dfont) \simeq \dfont$ et d'une action de $\Gamma$ semi-lin\'eaire continue et commutant \`a $\phi$. Rappelons que Fontaine a
construit dans \cite[A.3.4]{F90} un foncteur $T \mapsto \dfont(T)$ qui \`a toute $\OO_L$-repr\'esentation de $\gal$ associe un $(\phi,\Gamma)$-module
\'etale et que ce foncteur est une \'equivalence de cat\'egories. Si $V = L \otimes_{\OO_L} T$ est une repr\'esentation $p$-adique, alors on pose $\dfont(V) =  L \otimes_{\OO_L} \dfont(T)$. 

L'anneau $\calO$ est un $\phi(\calO)$-module libre de rang $p$,
dont une base est donn\'ee par $\{(1+X)^i\}_{0 \leq i \leq p-1}$. Si
$y \in \dfont$, alors on peut \'ecrire $y=\sum_{i=0}^{p-1} (1+X)^i
\phi(y_i)$ et on d\'efinit un op\'erateur $\psi : \dfont \to \dfont$ par la formule $\psi(y) = y_0$ si $y=\sum_{i=0}^{p-1} (1+X)^i \phi(y_i)$. 

Si $\dfont$ est un $(\phi,\Gamma)$-module \'etale sur $\calO$, alors Colmez
a d\'efini dans \cite[\S 4.5]{CL} un sous-$\OO_L[[X]]$-module $\ddiese$ de $\dfont$, qui est caract\'eris\'e par les propri\'et\'es suivantes : 
\begin{enumerate}
\item $\ddiese$ est un {\og treillis \fg} de $\dfont$ (voir \cite[\S 4.1]{CL});
\item quels que soient $x \in \dfont$ et $k \geq 0$, il existe $n(x,k) \geq 0$ tel que $\psi^n(x) \in \ddiese + p^k \dfont$ si $n \geq n(x,k)$;
\item l'op\'erateur $\psi$ induit une surjection de $\ddiese$ sur lui-m\^eme.
\end{enumerate}

Le calcul de $\ddiese(V)$ pour des repr\'esentations de dimension $1$ ne pose aucun probl\`eme et le calcul suivant sera utile dans la suite de cet article.

\begin{lemm}\label{calcidiot}
Si $W = k_L \cdot w$ est une $k_L$-repr\'esentation de $\gal$ de dimension $1$, donn\'ee par un caract\`ere $\omega^r \mu_y$, alors $\ddiese(W) = X^{-1} k_L[[X]] \cdot e$, avec $e=\alpha w$ o\`u $\alpha \in \Fpbar$ est tel que $\alpha^{p-1} = y$. 
\end{lemm}

Sous les hypoth\`eses du lemme ci-dessus, on pose alors $\dplus(W) = k_L[[X]] \cdot e$. Remarquons en passant que $\phi(e)=ye$ et que $\gamma(e)=\omega^r(\gamma)e$ si $\gamma \in \Gamma$. La notation $\dplus$ provient de la th\'eorie des repr\'esentations {\og de hauteur finie \fg} (qui n'intervient pas dans le reste de cet article).

\begin{lemm}\label{schur_pgmod}
Si $W$ est une $k_L$-repr\'esentation irr\'eductible de $\gal$ de dimension $\geq 2$, et si $M \subset \ddiese(W)$ est un sous-$k_L[[X]]$-module non-nul stable par  $\psi$, alors $M=\ddiese(W)$.

De m\^eme, si $W$ est une $k_L$-repr\'esentation de $\gal$ de dimension $1$, et si $M \subset \dplus(W)$ est un sous-$k_L[[X]]$-module non-nul stable par $\psi$, alors $M=\dplus(W)$.
\end{lemm}

\begin{proof}
Faisons tout d'abord le cas o\`u $W$ est de dimension $\geq 2$. Pour tout polyn\^ome $P \in k_L[X]$, on a $\dfont(W)^{P(\phi)=0} = 0$ parce que (voir la preuve du (iii) de la remarque 5.5 de \cite{CL}) on a : 
\[ \dfont(W)^{P(\phi)=0} \subset (\Fpbar
\otimes_{\Fp} W)^{\mathrm{Gal} (\Qpbar / \Qp(\mu_{p^\infty})) }
\subset \Fpbar \otimes_{\Fp} W^{ \mathrm{Gal} (\Qpbar / \Qp^{\mathrm{ab}}) } =  0. \]
Le lemme suit alors de la proposition 4.47 de \cite{CL} (ou plus
exactement de sa d\'emonstration, en remarquant que la d\'emonstration
n'utilise pas le fait que $\psi : M \to M$ est surjectif). 

Faisons maintenant le cas o\`u $W$ est de dimension $1$. Il s'agit de montrer qu'un id\'eal de $k_L[[X]]$ non-nul et stable par $\psi$ est \'egal \`a $k_L[[X]]$. Les id\'eaux de $k_L[[X]]$ sont de la forme $X^j k_L[[X]]$ pour $j \geq 0$ et un calcul facile montre que si $j \geq 1$, alors $\psi(X^j k_L[[X]])$ contient $X^{j-1} k_L[[X]]$ ce qui permet de conclure.
\end{proof}

\Subsection{Construction de repr\'esentations de $\B$}\label{constb}

Commen\c{c}ons ce paragraphe en rappelant que si $\Pi$ est une $k_L$-repr\'esentation compacte d'un groupe localement profini $G$, alors sa duale $\Pi^* = \injlim_{U \subset \Pi} \Hom(\Pi/U,k_L)$, o\`u $U$ parcourt les sous-$k_L$-espaces vectoriels ouverts de $\Pi$, est une $k_L$-repr\'esentation lisse de $G$ et que r\'eciproquement si $\Omega$ est une $k_L$-repr\'esentation lisse de $G$, alors sa duale $\Omega^* = \projlim_{H \subset G} \Hom(\Pi^H,k_L)$, o\`u $H$ parcourt les sous-groupes ouverts compacts de $G$, est une $k_L$-repr\'esentation compacte de $G$. De plus, on a $(\Pi^*)^* \simeq \Pi$ et $(\Omega^*)^* \simeq \Omega$ et les sous-espaces ferm\'es de $\Pi$ stables par $G$ correspondent aux quotients de $\Omega$ stables par $G$, et vice-versa.

Si $\dfont$ est un $(\phi,\Gamma)$-module, alors $\projlim_\psi \ddiese$ d\'enote l'ensemble des suites $v=(v_i)_{i \geq 0}$ telles que $\psi(v_{i+1})=v_i$ pour tout $i \geq 0$. Si $\dfont$ est un $(\phi,\Gamma)$-module sur $L$, on demande en plus que la suite $(v_i)_{i \geq 0}$ soit born\'ee pour la topologie $p$-adique.

On fixe un caract\`ere lisse $\chi$ de $\Qp^\times$ et on munit $\projlim_\psi \ddiese$ d'une action de $\B$ comme suit. Tout \'el\'ement $g \in \B$ peut s'\'ecrire comme produit :
\[ g = \begin{pmatrix} x & 0 \\ 0 & x \end{pmatrix} \cdot
\begin{pmatrix} 1 &  0 \\ 0 & p^j  \end{pmatrix} \cdot 
\begin{pmatrix} 1 &  0 \\ 0 & a  \end{pmatrix} \cdot 
\begin{pmatrix} 1 &  z \\ 0 & 1 \end{pmatrix}, \]
o\`u $x \in \Qp^\times$, $j\in\ZZ$, $a \in \Zp^{\times}$ et $z
\in \Qp$. Si $v=(v_i)_{i \geq 0} \in \projlim_{\psi} \ddiese$, alors on pose pour $i\geq 0$ :
\begin{align*}
\left( \begin{pmatrix} x & 0 \\ 0 & x \end{pmatrix} \star v \right)_i
& = \chi^{-1}(x) v_i;  \\
\left( \begin{pmatrix} 1 &  0 \\ 0 & p^j  
\end{pmatrix} \star v \right)_i & = v_{i-j} = \psi^j(v_i); \\
\left( \begin{pmatrix} 1 &  0 \\ 0 & a  
\end{pmatrix} \star v \right)_i & = \gamma_a^{-1}(v_i), 
\text{ o\`u $\gamma_a \in \Gamma$ est tel que $\eps(\gamma_a) = a$;}\\
\left( \begin{pmatrix} 1 &  z \\ 0 & 1  
\end{pmatrix} \star v \right)_i & = 
\psi^j((1+X)^{p^{i+j} z} v_{i+j}),\text{ pour $i+j \geq -{\rm val}(z)$.}
\end{align*}

\begin{defi}\label{defomega}
Si $W$ est une $k_L$-repr\'esentation irr\'eductible de $\gal$ de dimension $\geq 2$, 
et si $\chi : \Qp^\times \to k_L^\times$ est un caract\`ere lisse, alors on pose $\Omega_\chi(W) = (\projlim_\psi \ddiese(W))^*$ et si $W$ est une $k_L$-repr\'esentation de dimension $1$, alors on pose $\Omega_\chi(W) = (\projlim_\psi \dplus(W))^*$.
\end{defi}

\begin{prop}\label{omirred}
Si on munit $\Omega_\chi(W)$ d'une action de $\B$ en utilisant les formules ci-dessus, alors $\Omega_\chi(W)$ est une repr\'esentation lisse irr\'eductible dont le caract\`ere central est $\chi$.
\end{prop}

\begin{proof}
Comme $\projlim_\psi \ddiese(W)$ est compact, le fait que $\Omega_\chi(W)$ est une repr\'esen\-tation lisse r\'esulte des rappels que l'on a faits au d\'ebut du paragraphe. Le fait que le caract\`ere central est $\chi$ est \'evident. Enfin, pour montrer que $\Omega_\chi(W)$ est irr\'eductible, il suffit de montrer que $\projlim_\psi \ddiese(W)$ (ou bien $\projlim_\psi \dplus(W)$ si $\dim W = 1$) n'admet pas de sous-espace ferm\'e stable par $\B$ et non-trivial, ce que nous faisons maintenant (si $\dim W = 1$, il faut remplacer $\ddiese(W)$ par $\dplus(W)$ dans les calculs qui suivent). 

Soit $\pr_j : \projlim_\psi \ddiese(W) \to \ddiese(W)$ la projection $(v_n)_{n \geq 0} \mapsto v_j$. Si $M$ est un sous-espace ferm\'e et stable par $\B$ de $\projlim_\psi \ddiese(W)$, on note $M_j$ l'image de $\pr_j : M \to \ddiese(W)$. On voit que $M_j$ est un sous-$k_L[[X]]$-module non-nul de $\ddiese(W)$ stable par $\psi$ ce qui fait que, par le lemme \ref{schur_pgmod}, $M_j=\ddiese(W)$. On en d\'eduit que $M$ est dense dans $\projlim_\psi \ddiese(W)$ et donc finalement que $M = \projlim_\psi \ddiese(W)$ ce qui fait que $\projlim_\psi \ddiese(W)$ est bien topologiquement irr\'eductible et donc que $\Omega_\chi(W)$ est irr\'eductible. 

\end{proof}

\begin{prop}\label{schur}
Si $W_1$ et $W_2$ sont deux $k_L$-repr\'esentations irr\'eductibles de $\gal$ et si $f : \Omega_{\chi_1}(W_1) \to \Omega_{\chi_2}(W_2)$ est une application $\B$-\'equivariante et non-nulle, alors $\chi_1 = \chi_2$ et $W_1 \simeq W_2$ et $f$ est scalaire. 
\end{prop}

\begin{proof}
Il est imm\'ediat que $\chi_1 = \chi_2$ et on se ram\`ene donc \`a montrer que si $W_1$ et $W_2$ sont deux repr\'esentations irr\'eductibles telles qu'il existe une application $\left(\begin{smallmatrix} 1 & \ast \\ 0 & \ast \end{smallmatrix}\right)$-\'equivariante continue et non-nulle $f : \projlim_\psi \ddiese(W_1) \to \projlim_\psi \ddiese(W_2)$, alors $W_1 \simeq W_2$. 

La d\'emonstration est analogue \`a celle de la proposition 3.4.3 de
\cite{BB} et nous la faisons pour $\dim W_i \geq 2$; quand $\dim W_i = 1$, il suffit de remplacer $\ddiese(W_i)$ par $\dplus(W_i)$ dans les calculs qui suivent.

Notons comme ci-dessus $\pr_0 : \projlim_\psi \ddiese(W)
\to \ddiese(W)$ la projection $(v_n)_{n \geq 0} \mapsto v_0$.
Commen\c{c}ons par montrer que si $v=(v_n)_{n \geq 0} \in \projlim_\psi
\ddiese(W_1)$, alors $\pr_0 \circ f(v)$ ne d\'epend que de $v_0 = \pr_0(v)$. Pour cela, soit $K_n$ l'ensemble des $v \in \projlim_\psi \ddiese(W_1)$ dont les $n$ premiers termes sont nuls, ce qui fait que pour $n \geq 1$, $K_n$ est un sous-$k_L[[X]]$-module ferm\'e et stable par $\psi$ de 
$\projlim_\psi \ddiese(W_1)$ et que $\psi(K_n)=K_{n+1}$. On en d\'eduit que $\pr_0 \circ f (K_n)$ est un sous-$k_L[[X]]$-module ferm\'e et stable par $\psi$ de $\ddiese(W_2)$. Le lemme \ref{schur_pgmod} implique alors que soit $\pr_0 \circ f (K_n) =
0$, soit $\pr_0 \circ f (K_n) = \ddiese(W_2)$. Enfin, on voit que
$\psi(\pr_0 \circ f (K_n)) = \pr_0 \circ f (K_{n+1})$ d'une part et que $\pr_0 \circ f (K_n) = 0$ si $n \gg 0$ par continuit\'e d'autre part. Cela implique que $\pr_0 \circ f (K_n) = 0$ pour tout $n \geq 1$ et donc que si $v_0 = 0$, alors $\pr_0 \circ f(v) = 0$. 

Pour tout $w \in \ddiese(W_1)$, soit $\widetilde{w}$ un \'el\'ement de 
$\projlim_{\psi} \ddiese(W_1)$ tel que $\widetilde{w}_0 = w$. Les calculs pr\'ec\'edents montrent que l'application $h :  \ddiese(W_1) \to \ddiese(W_2)$ donn\'ee par $h(w) = \pr_0 \circ f (\widetilde{w})$ est bien d\'efinie, et
qu'elle est $k_L[[X]]$-lin\'eaire et commute \`a $\psi$ et \`a l'action
de $\Gamma$. Par les propositions 4.7 et 4.55 de \cite{CL}, elle s'\'etend en une application de  $(\phi,\Gamma)$-modules $h: \dfont(W_1) \to \dfont(W_2)$ et par fonctorialit\'e, on en d\'eduit qu'il existe une application non-nulle et $\gal$-\'equivariante de $W_1$ dans $W_2$, ce qui fait que par le lemme de Schur, on a $W_1 \simeq W_2$ et que $f$ est scalaire.
\end{proof}

Si $W$ est une $k_L$-repr\'esentation de dimension $1$, le lemme \ref{calcidiot} montre que $\ddiese(W) = X^{-1} k_L[[X]] \cdot e$. On appelle $\res$ l'application de $\ddiese(W)$ dans $k_L$ qui \`a $y = \sum_{j=-1}^\infty f_j X^j \cdot e$ associe $f_{-1}$ et on note toujours $\res$ l'application de $\projlim_\psi \ddiese(W)$ dans $k_L$ qui \`a $(y_0,y_1,\hdots)$ associe $\res(y_0)$. 

Si $\chi_1$ et $\chi_2$ sont deux caract\`eres de $\Qp^\times$, on note $\chi_1 \otimes \chi_2$ le caract\`ere de $\B$ d\'efini par la formule :
\[ \chi_1 \otimes \chi_2 : \begin{pmatrix} a & b \\ 0 & d \end{pmatrix} \mapsto \chi_1(a)\chi_2(d). \]

\begin{prop}\label{plusetdiese}
Si $\eta_W$ d\'enote le caract\`ere de $\Qp^\times$ associ\'e \`a $W$ par le corps de classes, alors l'application $\res$ induit par dualit\'e une suite exacte de $\B$-repr\'esen\-tations :
\[ 0 \to \chi \omega \eta_W^{-1} \otimes \omega^{-1} \eta_W \to (\projlim_\psi \ddiese(W))^* \to \Omega_\chi(W) \to 0. \]
\end{prop}

\begin{proof}
Il s'agit de montrer que si $W = \omega^r \mu_\lambda$, alors on a une suite exacte :
\[ 0 \to \projlim_\psi \dplus(W) \to \projlim_\psi \ddiese(W) \to
\chi^{-1} \omega^{r-1} \mu_\lambda \otimes \omega^{1-r} \mu_{\lambda^{-1}} \to 0. \]

Si l'on pose $y=(\lambda^n X^{-1} e)_{n \geq 0}$, alors $y \in \projlim_\psi \ddiese(W)$ et $\res(y)=1$ ce qui fait qu'il suffit de calculer $\res(b \star y)$ pour $b \in \B$. On trouve que :
\begin{align*}
\res \left(  \left(\begin{smallmatrix}  x & 0 \\ 0 & x \end{smallmatrix}\right)  \star y \right) & = \chi^{-1}(x), \\
\res \left(  \left(\begin{smallmatrix} 1 & 0 \\ 0 & p \end{smallmatrix}\right)  \star y \right) & = \res (( \lambda^{n-1}X^{-1} e )_{n \geq 0}) = \lambda^{-1},
\end{align*}

\begin{align*}
\res \left(  \left(\begin{smallmatrix} 1 & 0 \\ 0 & a \end{smallmatrix}\right)  \star y \right) & = \res ((\omega^{r-1}(\gamma_a^{-1}) (\lambda^n +\mathrm{O}(X)) X^{-1} e)_{n \geq 0})  = a^{1-r}, \\
\res \left(  \left(\begin{smallmatrix}  1 & z \\ 0 & 1 \end{smallmatrix}\right)  \star y \right) & = 1,
\end{align*}
ce qui permet de conclure.
\end{proof}

\begin{rema}
Nous allons voir dans la suite de cet article que si $\dim W \in \{ 1,2 \}$, alors les repr\'esentations $\Omega_\chi(W)$ sont les restrictions \`a $\B$ de repr\'esentations de $\G$. Quand $\dim W \geq 3$, ce n'est plus le cas et on obtient des repr\'esentations lisses irr\'eductibles propres \`a $\B$. Obtient-on ainsi toutes les repr\'esentations lisses irr\'eductibles de $\B$?
\end{rema}

\Subsection{Repr\'esentations de $\G$}\label{repg}

L'objet de ce paragraphe est de rappeler bri\`evement certaines des constructions de Breuil et Colmez. Il s'agit d'associer une repr\'esentation $\Pi(V)$ de $\G$ \`a une repr\'esentation $p$-adique $V$ irr\'eductible de dimension $2$. Cette repr\'esentation $\Pi(V)$ a \'et\'e conjecturalement construite pour les repr\'esentations cristallines et semi-stables par Breuil (voir \cite{Br2} pour le cas cristallin et \cite{BRst} pour le cas semi-stable). Le fait que la repr\'esentation $\Pi(V)$ est non-nulle, irr\'eductible et admissible (au sens de \cite{ST}) a \'et\'e d\'emontr\'e par Colmez pour $V$ semi-stable (voir \cite{CL}; des cas particuliers se trouvent dans \cite{BRst}) puis par Breuil et l'auteur pour $V$ cristalline (voir \cite{BB}; des cas particuliers se trouvent dans \cite{Br2}). Enfin, ces r\'esultats ont \'et\'e \'etendus aux repr\'esentations triangulines par Colmez (voir \cite{CT}). Le th\'eor\`eme ci-dessous rassemble les r\'esulats dont nous nous servirons\footnote{Ces r\'esultats ne sont pas r\'edig\'es dans les cas {\og exceptionnels \fg} ($V$ cristalline non $\phi$-semi-simple ou bien $V$ non-cristalline mais le devenant sur une extension ab\'elienne).}.

\begin{theo}\label{colmod}
Si $V$ est une repr\'esentation trianguline irr\'eductible de dimension $2$, alors :
\begin{enumerate}
\item on a $\Pi(V)^* \simeq_B \projlim_\psi \ddiese(V)$;
\item si $T$ est un r\'eseau de $V$, alors $(\projlim_\psi \ddiese(T))^*$ est un r\'eseau de $\Pi(V)$ stable sous l'action de $\B$ et commensurable aux r\'eseaux de $\Pi(V)$;
\item on a $\overline{\Pi}(V)^* \simB \projlim_\psi \ddiese(\overline{V})$.
\end{enumerate}
\end{theo}

\begin{proof}
Le (1) est d\'emontr\'e dans \cite{CL} pour une repr\'esentation semi-stable, dans \cite{BB} pour une repr\'esentation cristalline et enfin dans \cite{CT} dans le cas d'une repr\'esentation trianguline\footnote{Il faut tout de m\^eme faire attention au fait que les normalisations changent d'un article \`a l'autre et il convient donc \'eventuellement de tordre la formule par un caract\`ere.}. Le (2) est un r\'esultat facile (voir \cite[\S 5.4]{BB} par exemple) et le (3) suit du (1) et du (2) ainsi que de la proposition 4.50 de \cite{CL} et du principe de Brauer-Nesbitt (ce qui explique qu'il faut semi-simplifier). 
\end{proof}

Rappelons \`a pr\'esent la classification des $k_L$-repr\'esentations lisses irr\'eductibles de $\G$ qui admettent un caract\`ere central. Si $r \in \{0,\hdots,p-1\}$, si $\lambda \in \Fpbar$ et si $\chi : \Qp^\times \to k_L^{\times}$ est un caract\`ere continu, on pose : 
\[ \pi(r,\lambda,\chi) = \left( \frac{\ind_{\K\Qp^\times}^{\G}
    \mathrm{Sym}^r k_L^2}{T-\lambda} \right) \otimes (\chi \circ
\det), \] 
o\`u $T$ est un certain op\'erateur de Hecke (voir \cite{BL}). 

Si $(r,\lambda) \notin \{ (0,\pm 1); (p-1,\pm 1) \}$, alors $\pi(r,\lambda,\chi)$ est irr\'eductible (pour $\lambda \neq 0$, c'est un r\'esultat de \cite{BL} et pour $\lambda = 0$, c'est un r\'esultat de \cite{Br1}). Pour $\lambda=0$, les entrelacements entre les $\pi(r,0,\chi)$ sont les suivants : 
\[ \pi(r,0,\chi) \simeq \pi(r,0,\chi \mu_{-1}) \simeq \pi(p-1-r,0,\chi \omega^r) \simeq \pi(p-1-r,0,\chi \omega^r \mu_{-1}). \]

Si $(r,\lambda) \in \{ (0,\pm 1); (p-1,\pm 1) \}$, alors la semi-simplifi\'ee de $\pi(r,\lambda,\chi)$ est somme d'un caract\`ere et de la tordue d'une repr\'esentation appel\'ee la sp\'eciale (qui correpond \`a $r=0$, $\lambda=1$ et $\chi=1$). Enfin, toute repr\'esentation lisse irr\'eductible de $\G$ qui admet un caract\`ere central\footnote{Dans la suite, on suppose tout le temps sans le dire que les repr\'esentations lisses irr\'eductibles de $\G$ que l'on consid\`ere admettent un caract\`ere central.} est dans la liste ci-dessous :
\begin{enumerate}
\item les caract\`eres $\chi \circ \det$;
\item la sp\'eciale tordue par un caract\`ere $\mathrm{Sp} \otimes (\chi \circ \det)$;
\item les s\'eries principales $\pi(r,\lambda,\chi)$ avec $\lambda \neq 0$ et $(r,\lambda) \notin \{ (0,\pm 1); (p-1,\pm 1) \}$;
\item les supersinguli\`eres $\pi(r,0,\chi)$.
\end{enumerate}

Nous rappelons maintenant le calcul de $\overline{\Pi}(V)$ pour certaines repr\'esentations cristallines. Si $k \geq 2$ et si $a_p \in \MM_L$, alors on d\'efinit un $\phi$-module filtr\'e $D_{k,a_p}$ par $D_{k,a_p} = L e \oplus L
f$ o\`u : 
\[ \begin{cases} \phi(e) = p^{k-1} f \\
\phi(f) = -e + a_p f 
\end{cases}
\text{et}\quad
\Fil^i D_{k, a_p} = \begin{cases}
D_{k, a_p} & \text{si $i \leq 0$,} \\
L e & \text{si $1 \leq i \leq k-1$,} \\
0 & \text{si $i \geq k$.}
\end{cases} \]
Ce $\phi$-module filtr\'e est admissible, et  
on sait (par le th\'eor\`eme principal de \cite{CF})
qu'il existe alors une repr\'esentation cristalline
$V_{k,a_p}$ de $\gal$
telle que $\dcris(V^*_{k,a_p}) = D_{k,a_p}$ (on passe au
dual pour que les notations soient compatibles avec celles de
\cite{Br2} et \cite{BLZ}; remarquons tout de m\^eme que 
l'on a $V_{k,a_p}^* = V_{k,a_p}(1-k)$). 

Soit $\Pi_{k,a_p}$ la repr\'esentation localement alg\'ebrique construite par Breuil dans \cite{Br2} :
\[ \Pi_{k,a_p} = \frac{\ind_{\K\Qp^\times}^{\G} \mathrm{Sym}^{k-2} L^2}{T-a_p}, \]
o\`u $T$ est un certain op\'erateur de Hecke (voir \cite{Br2}).
 
Il est montr\'e dans \cite{BB} que $\Pi_{k,a_p}$ admet un $\OO_L$-r\'eseau de type fini sur $\OO_L[\G]$ et que sa compl\'etion est isomorphe \`a $\Pi(V_{k,a_p})$. Breuil a d'autre part d\'etermin\'e $\overline{\Pi}_{k,a_p}$ pour $k \leq 2p+2$ (voir le paragraphe \ref{calc}) et nous aurons besoin dans la suite du cas particulier ci-dessous (c'est le corollaire 5.1 de \cite{Br2}).

\begin{prop}\label{breuilmod}
Si $k \leq p+1$, alors $\overline{\Pi}_{k,0} \simeq \pi(k-2,0,1)$.
\end{prop}

\section{Etude de certaines repr\'esentations de $\B$}\label{repb}

L'objet de ce chapitre est de faire le lien entre les diff\'erentes constructions de repr\'esentations lisses irr\'eductibles de $\B$ (soit comme restriction de repr\'esentations lisses irr\'eductibles de $\G$, soit directement comme $\Omega_\chi(W)$). 

\Subsection{Induites paraboliques}\label{indpar}

On commence par faire l'\'etude des induites paraboliques. Si $\chi_1$ et $\chi_2$ sont deux caract\`eres de $\Qp^\times$, alors l'induite parabolique $\Indbg (\chi_1 \otimes \chi_2)$ est par d\'efinition l'ensemble des fonctions localement constantes $\sigma : \G \to k_L$ qui v\'erifient $\sigma(bg)  = (\chi_1 \otimes \chi_2)(b)\sigma(g)$. 

Si $\sigma \in \Indbg (\chi_1 \otimes \chi_2)$, alors l'application $\sigma \mapsto \sigma(\Id)$ nous donne une application de $\Indbg (\chi_1 \otimes \chi_2)$ dans $\chi_1 \otimes \chi_2$ dont on note  $\Indbg (\chi_1 \otimes \chi_2)_0$ le noyau ce qui fait que l'on a une suite exacte de $\B$-repr\'esentations :
\[ 0 \to \Indbg (\chi_1 \otimes \chi_2)_0 \to \Indbg (\chi_1 \otimes \chi_2) \to \chi_1 \otimes \chi_2 \to 0. \]
 
Soit $\LC_0(\Qp,k_L)$ l'ensemble des fonctions $f : \Qp \to k_L$
localement constantes \`a support compact. Si $\sigma \in \Indbg (\chi_1 \otimes \chi_2)$, alors on associe \`a $\sigma$ une fonction $f_\sigma \in \LC_0(\Qp,k_L)$ via la recette :  
\[ f_\sigma(x) = \sigma \left(  \begin{pmatrix} 0 & 1 \\ -1 & x \end{pmatrix} \right). \]
La d\'ecomposition de Bruhat montre que l'application $\sigma \mapsto f_\sigma$ est une bijection de $\Indbg (\chi_1 \otimes \chi_2)_0$ dans $\LC_0(\Qp,k_L)$ et on a alors :
\begin{align*}
f \left( \begin{pmatrix} a & b \\ 0 & d \end{pmatrix} \star \sigma \right) (x) 
& = \sigma \left( \begin{pmatrix} 0 & 1 \\ -1 & x \end{pmatrix} 
\begin{pmatrix} a & b \\ 0 & d \end{pmatrix} \right) \\ \tag{eq1}
& = \sigma \left( \begin{pmatrix} d & 0 \\ 0 & a \end{pmatrix} 
\begin{pmatrix} 0 & 1 \\ -1 &  \frac{dx-b}{a} \end{pmatrix} \right) \\
& = \chi_1(d)\chi_2(a) f_\sigma \left( \frac{dx-b}{a} \right). 
\end{align*}

\begin{theo}\label{omdim1}
Si $\eta_W$ d\'enote le caract\`ere de $\Qp^\times$ associ\'e \`a $W$ par le corps de classes, alors $\Omega_\chi(W) \simeq \Indbg (\eta_W \otimes \chi \eta_W^{-1})_0$.
\end{theo}

\begin{proof}
Ecrivons $W = \omega^r \mu_\lambda$, ce qui fait que l'on doit montrer que :
\[ \Omega_\chi(W)\simeq \Indbg (\omega^r \mu_\lambda \otimes \chi \omega^{-r} \mu_{\lambda^{-1}})_0.  \] 
Nous allons d\'efinir une application : \[ \projlim_\psi \dplus(W) \to \Indbg (\omega^r \mu_\lambda \otimes \chi \omega^{-r} \mu_{\lambda^{-1}})_0^*,  \] 
et montrer que c'est un isomorphisme $\B$-\'equivariant.

Rappelons tout d'abord la d\'efinition de la transform\'ee d'Amice (tout au moins la version dont nous avons besoin). Si $\nu$ est une forme lin\'eaire $\nu :  \LC(\Zp,k_L) \to k_L$ (on dira que $\nu$ est une mesure sur $\Zp$), alors sa transform\'ee d'Amice $\calA(\nu) \in k_L[[X]]$ est d\'efinie par : 
\[ \calA(\nu)(X) = \sum_{n=0}^\infty \nu\left( z \mapsto \binom{z}{n} \right) X^n = \nu(z \mapsto (1+X)^z), \]
et l'application $\calA :  \LC(\Zp,k_L)^* \to k_L[[X]]$ est alors un isomorphisme.

Fixons une base $e$ de $\dplus(W)$ telle que $\phi(e) = \lambda e$ et $\gamma(e) = \omega^r(\gamma)$ (voir le lemme \ref{calcidiot}); si $y \in \projlim_\psi \dplus(W)$, alors $y=(y_i)_{i \geq 0}$ et on a $y_i = f_i e$ avec $f_i \in k_L[[X]]$. On voit que $\psi(\lambda^{-i} f_i) =  \lambda^{-(i-1)} f_{i-1}$ et on d\'efinit une mesure $\nu_{y,i}$ sur $\Zp$ en demandant que $\calA(\nu_{y,i}) = \lambda^{-i} f_i$ ce qui nous permet de d\'efinir une mesure $\nu_y$ sur $\Qp$ en imposant que si $f \in \LC_0(\Qp,k_L)$ est \`a support dans $p^{-i} \Zp$, alors :
\[ \int_{\Qp} f d\nu_y = \int_{\Zp} f(p^{-i} z) d\nu_{y,i}. \]
Le terme de droite ne d\'epend pas de $i \gg 0$, en vertu de la formule :
\[ \int_{\Zp} f(z) d\psi(\nu) = \int_{p\Zp} f(p^{-1}z) d\nu, \]
et $\int_{\Qp} f d\nu_y$ est donc bien d\'efinie. Il est clair que l'application $y \mapsto \nu_y$ est une bijection qui nous permet d'identifier, via la recette $\sigma \mapsto f_\sigma$ rappel\'ee au d\'ebut de ce paragraphe, $\Omega_\chi(W)$ avec le dual de l'induite parabolique.

Pour terminer la d\'emonstration du th\'eor\`eme, nous devons montrer que cette identification est $\B$-\'equivariante. L'action de $\B$ sur les mesures est donn\'ee par l'action correspondante sur leurs transform\'ees d'Amice et on trouve les formules suivantes : 
\begin{align*}
\int_{\Qp} f(z) d\nu_{\left(\begin{smallmatrix}  x & 0 \\ 0 & x \end{smallmatrix}\right) \star y} & = 
\chi^{-1}(x) \int_{\Qp} f(z) d\nu_y \\
\int_{\Qp} f(z) d\nu_{\left(\begin{smallmatrix}  1 & 0 \\ 0 & p \end{smallmatrix}\right) \star y} & = 
\lambda^{-1} \int_{\Qp} f(p^{-1}z) d\nu_y \\
\int_{\Qp} f(z) d\nu_{\left(\begin{smallmatrix}  1 & 0 \\ 0 & d \end{smallmatrix}\right) \star y} & = 
d^{-r} \int_{\Qp} f(d^{-1}z) d\nu_y \quad(\text{avec $d \in \Zp^\times$}) \\
\int_{\Qp} f(z) d\nu_{\left(\begin{smallmatrix}  1 & b \\ 0 & 1 \end{smallmatrix}\right) \star y} & = 
\int_{\Qp} f(z+b) d\nu_y.
\end{align*}
En comparant ces formules avec (eq1), on trouve bien que : 
\[ \int_{\Qp} f_{g \star \sigma} (z) d\nu_{g \star y} = \int_{\Qp} f_{\sigma}(z) d\nu_y \]
si $g \in \B$, $y \in \Omega_\chi(W)$ et $\sigma \in \Indbg (\omega^r \mu_\lambda \otimes \chi \omega^{-r} \mu_{\lambda^{-1}})_0$.
\end{proof}

Pour terminer, rappelons le lien entre les induites paraboliques et les $\pi(r,\lambda,\chi)$ (c'est le th\'eor\`eme 30 de \cite{BL2}).

\begin{prop}\label{pietind}
On a $\pi(r,\lambda,\chi) \simG \Indbg(\chi \mu_{\lambda^{-1}} \otimes \chi \mu_{\lambda} \omega^r)$. 
\end{prop}

\begin{rema}\label{spe}
En particulier, la sp\'eciale $\mathrm{Sp}$ est isomorphe \`a $\Indbg(1 \otimes 1)_0$. Dans ce cas, la suite exacte de repr\'esentations de $\B$ :
\[ 0 \to \Indbg(1 \otimes 1)_0 \to \Indbg(1 \otimes 1) \to 1 \otimes 1 \to 0 \] 
est d'ailleurs scind\'ee (le $k_L$-espace vectoriel engendr\'e par la fonction $g \mapsto 1$ est un suppl\'ementaire $\B$-stable de $\mathrm{Sp}$). Il est amusant de remarquer que la suite exacte de repr\'esentations de $\G$ :
\[ 0 \to \mathrm{Id} \to \Indbg(1 \otimes 1) \to \mathrm{Sp} \to 0 \] 
n'est \emph{pas} scind\'ee (voir la d\'emonstration du th\'eor\`eme 29 de \cite{BL}). 
\end{rema}

\Subsection{Supersinguli\`eres}\label{supsing}

Dans ce court paragraphe, nous d\'emontrons l'analogue du th\'eor\`eme \ref{omdim1} pour les repr\'esentations de dimension $2$ mais contrairement \`a la d\'emonstration du th\'eor\`eme \ref{omdim1}, l'argument n'est pas direct et repose de mani\`ere essentielle sur des calculs en caract\'eristique $0$. Il serait int\'eressant de disposer d'une d\'emonstration purement en caract\'eristique $p$.

\begin{theo}\label{omdim2}
Si $W \simeq \rho(r,\chi)$, alors $\Omega_{\omega^r \chi^2}(W)$ est isomorphe \`a la restriction au Borel de $\pi(r,0,\chi)$.  
\end{theo}

\begin{proof}
Un calcul facile montre que $\overline{V}_{r+2,0} \simeq \rho(r,1)$ et le (3) du th\'eor\`eme \ref{colmod} montre alors que $\projlim_\psi \ddiese(\rho(r,\chi)) \simB \overline{\Pi}(V_{r+2,0} \otimes \chi)^*$ ce qui fait que, comme le caract\`ere central de $\overline{\Pi}_{r+2,0} \otimes (\chi \circ \det)$ est $\omega^r \chi^2$, 
on a : 
\[ \Omega_{\omega^r \chi^2}(\rho(r,\chi)) \simB \overline{\Pi}_{r+2,0} \otimes (\chi \circ \det) \simeq \pi(r,0,\chi), \] 
par la proposition \ref{breuilmod}.   
\end{proof}

\Subsection{Applications aux repr\'esentations de $\G$}\label{appli}

Pour terminer ce chapitre, nous donnons deux applications des identifications pr\'ec\'e\-dentes \`a l'\'etude des repr\'esentations de $\B$ et de $\G$.

\begin{theo}\label{resbor}
Si $\Pi$ est une $k_L$-repr\'esentation lisse irr\'eductible de $\G$ admettant un caract\`ere central, alors :
\begin{enumerate}
\item si $\Pi$ est un caract\`ere, ou la sp\'eciale, ou supersinguli\`ere, alors sa restriction \`a $\B$ est toujours irr\'eductible;
\item si $\Pi$ est une s\'erie principale, alors sa restriction \`a $\B$ est une extension du caract\`ere induisant $\chi_1 \otimes \chi_2$ par la repr\'esentation irr\'eductible $\Indbg(\chi_1 \otimes \chi_2)_0$.
\end{enumerate} 
\end{theo}

\begin{proof}
Montrons tout d'abord le (1). Comme la sp\'eciale est isomorphe \`a $\Indbg(1 \otimes 1)_0$ on renvoie au (2) et le cas des caract\`eres est trivial, ce qui nous laisse les supersinguli\`eres. Comme $\pi(r,0,\chi) \simeq \Omega_{\omega^r \chi^2}(\rho(r,\chi))$, le fait que cette repr\'esentation est irr\'eductible en restriction au Borel suit de la proposition \ref{omirred}.

Passons \`a pr\'esent au (2). On a vu au paragraphe \ref{indpar} que l'on a une suite exacte : 
\[ 0 \to \Indbg (\chi_1 \otimes \chi_2)_0 \to \Indbg (\chi_1 \otimes \chi_2) \to \chi_1 \otimes \chi_2 \to 0, \]
ce qui montre l'assertion concernant l'extension, et comme $\Indbg (\chi_1 \otimes \chi_2)_0 \simeq \Omega_{\chi_1 \chi_2}(\chi_1)$ par le th\'eor\`eme \ref{omdim1}, le fait que cette repr\'esentation est irr\'eductible suit encore de la proposition \ref{omirred}.
\end{proof}

\begin{theo}\label{sameresbor}
Si $\Pi_1$ et $\Pi_2$ sont deux $k_L$-repr\'esentations lisses semi-simples et de longueur finie de $\G$ (dont les composantes irr\'eductibles admettent un caract\`ere central) et dont les semi-simplifications des restrictions \`a $\B$ sont isomorphes, alors $\Pi_1$ et $\Pi_2$ sont d\'ej\`a isomorphes en tant que repr\'esentations de $\G$.
\end{theo}

\begin{proof}
Une repr\'esentation semi-simple et de longueur finie de $\G$ est une somme directe de caract\`eres $\chi \circ \det$, de tordues de la sp\'eciale, de s\'eries principales $\Indbg (\chi_1 \otimes \chi_2)$ avec $\chi_1 \neq \chi_2$ et de supersinguli\`eres. 

La restriction \`a $\B$ du caract\`ere $\chi \circ \det$ est $\chi \otimes \chi$ et par le th\'eor\`eme \ref{resbor}, celle de $\Indbg (\chi_1 \otimes \chi_2)$ est une extension de $\chi_1 \otimes \chi_2$ (avec n\'ecessairement $\chi_1 \neq \chi_2$) par une repr\'esentation de dimension infinie alors que les restrictions des tordues de la sp\'eciale et des supersinguli\`eres sont irr\'eductibles et de dimension infinie. On voit donc que si l'on a deux repr\'esentations $\Pi_1$ et $\Pi_2$ qui satisfont les conditions du th\'eor\`eme, alors elles contiennent les m\^emes caract\`eres et s\'eries principales. 

Pour terminer, il faut voir que les tordues des sp\'eciales et les supersinguli\`eres sont d\'etermin\'ees par leur restriction au Borel. Elles sont toutes de la forme $\Omega_\chi(W)$ (avec $W$ de dimension $1$ pour la sp\'eciale par le th\'eor\`eme \ref{omdim1} et $W$ de dimension $2$ pour les supersinguli\`eres par le th\'eor\`eme \ref{omdim2}) et la proposition \ref{schur} montre qu'il n'y a pas de $\B$-entrelacements non-triviaux entre deux telles repr\'esentations de $\G$ diff\'erentes.
\end{proof}

\section{Correspondances continues et modulo $p$}\label{main}

Dans ce chapitre, nous d\'emontrons le r\'esultat principal de cet article (la correspondance entre r\'eduction modulo $p$ des repr\'esentations galoisiennes $V$ et des repr\'esentations $\Pi(V)$ de $\G$) et ensuite nous rappelons les cas o\`u le calcul de cette r\'eduction a \'et\'e fait.

\Subsection{D\'emonstration de la conjecture de Breuil}\label{demo}

Dans ce paragraphe, nous d\'emontrons le r\'esultat principal de cet article. Le th\'eor\`eme ci-dessous est (pour $V$ cristalline) la conjecture 1.2 de \cite{Br2} (attention au fait que $\mu_\lambda = \mathrm{nr}(\lambda^{-1})$).

\begin{theo}\label{youpla}
Si $V$ est une repr\'esentation trianguline irr\'educible, alors :
\begin{align*}
\overline{V} = \rho(r,\chi) & \Leftrightarrow \overline{\Pi}(V) = \pi(r,0,\chi) \\
\overline{V} = \begin{pmatrix} \mu_\lambda \omega^{r+1} & 0 \\ 0 & \mu_{\lambda^{-1}} \end{pmatrix} \otimes \chi & \Leftrightarrow \overline{\Pi}(V) \simG \pi(r,\lambda,\chi) \oplus \pi([p-3-r], \lambda^{-1}, \omega^{r+1} \chi)
\end{align*}
\end{theo}

Ici $[p-3-r]$ est l'unique entier appartenant \`a $\{ 0, \hdots, p-2 \}$ et qui est congruent \`a $r$ modulo $p-1$.

\begin{proof}
Le (3) du th\'eor\`eme \ref{colmod} nous dit que l'on a :
\[ \tag{eq2} \overline{\Pi}(V)^* \simB \projlim_\psi \ddiese(\overline{V}). \] 

Supposons tout d'abord que l'on connaisse $\overline{V}$. Le th\'eor\`eme \ref{omdim1} et la proposition \ref{plusetdiese} quand $\overline{V}$ est r\'eductible, ou bien le th\'eor\`eme \ref{omdim2} quand $\overline{V}$ est irr\'eductible, ainsi que la formule (eq2), montrent que $\overline{\Pi}(V)$ est une repr\'esentation de $\G$ dont la restriction \`a $\B$ co\"{\i}ncide avec celle qui est pr\'edite par le th\'eor\`eme, et le th\'eor\`eme \ref{sameresbor} permet alors de conclure.

Supposons maintenant que l'on connaisse $\overline{\Pi}(V)$; on peut conclure comme ci-dessus (nous laissons cela en exercice) mais il y a plus explicite. Si $\overline{\Pi}(V)$ est une repr\'esentation supersinguli\`ere, alors par le th\'eor\`eme \ref{resbor} sa restriction \`a $\B$ est toujours irr\'eductible ce qui fait que (\`a cause de la formule (eq2)) $\overline{V}$ est irr\'eductible et le th\'eor\`eme \ref{omdim2} et la proposition \ref{schur} permettent de conclure. Si $\overline{\Pi}(V)$ est du deuxi\`eme type, alors $\overline{V}$ est forc\'ement r\'eductible et la formule (eq2) ainsi que la proposition \ref{plusetdiese} montrent que l'on peut lire $\overline{V}$ sur la restriction \`a $\B$ de $\overline{\Pi}(V)$ (c'est la tordue par $\omega$ de la partie de dimension finie de cette restriction, c'est cela qui est explicite).
\end{proof}

\begin{rema}\label{simpler}
La d\'emonstration ci-dessus montre qu'en fait, la restriction \`a $\B$ de $\overline{\Pi}(V)$ suffit \`a d\'eterminer $\overline{V}$. Dans les applications, c'est plut\^ot $\overline{\Pi}(V)$ que l'on calcule, et cela a pour effet de simplifier ces calculs\footnote{Nous reviendrons l\`a-dessus dans un prochain travail.}.
\end{rema}

\begin{rema}\label{ghost}
Si $W$ est de dimension $1$, on voit que la repr\'esentation $\Omega_\chi(W)$ est d'une part un sous-objet de $\Indbg(\eta_W \otimes \chi \eta_W^{-1})$ et d'autre part un quotient de $(\projlim_\psi \ddiese(W))^*$. On pourrait donc croire que $\Indbg(\eta_W \otimes \chi \eta_W^{-1}) \simB (\projlim_\psi \ddiese(W))^*$ mais en fait on voit que : 
\begin{align*}
\Indbg(\eta_W \otimes \chi \eta_W^{-1}) & \simB  \Omega_\chi(W) \oplus  (\eta_W \otimes \chi \eta_W^{-1}) \\
(\projlim_\psi \ddiese(W))^* & \simB \Omega_\chi(W) \oplus  (\chi \omega \eta_W^{-1} \otimes  \omega^{-1} \eta_W)
\end{align*}
et dans la correspondance entre $\mu_\lambda \omega^{r+1}\chi \oplus \mu_{\lambda^{-1}} \chi$ et $\pi(r,\lambda,\chi) \oplus \pi([p-3-r], \lambda^{-1}, \omega^{r+1} \chi)$, les termes se {\og croisent \fg} et il y a donc une sorte d'{\og entrelacement 
fant\^ome \fg} entre $\pi(r,\lambda,\chi)$ et $\pi([p-3-r], \lambda^{-1}, \omega^{r+1} \chi)$ (et il y a d'ailleurs un v\'eritable entrelacement en $\ell$-adique pour $\ell \neq p$, voir la remarque 4.2.5 de \cite{Br1}).
\end{rema}

\Subsection{R\'eduction modulo $p$ de quelques repr\'esentations galoisiennes}\label{calc}

L'objet de ce paragraphe est de rappeler les r\'esultats des calculs qui ont \'et\'e faits par diff\'erents auteurs. Il s'agit soit du calcul de $\overline{V}$, soit du calcul de $\overline{\Pi}(V)$, pour certaines repr\'esentations $V$. 

Commen\c{c}ons par les repr\'esentations cristallines $V_{k,a_p}$ d\'efinies \`a la fin du paragraphe \ref{repg}. Le calcul de $\overline{V}_{k,a_p}$ remonte \`a Fontaine et Serre (en raison du lien avec les formes modulaires); la th\'eorie de Fontaine-Laffaille (voir \cite{FL82}) permet de faire ce calcul pour tout $k \leq p$. Quand $k \leq p+1$, c'est un r\'esultat d'Edixhoven (sous l'hypoth\`ese que la repr\'esentation provienne d'une forme modulaire). Enfin, quand $k \geq p+1$, la th\'eorie des modules de Wach permet de calculer $\overline{V}_{k,a_p}$ si $\val(a_p)$ est assez grand (voir \cite{BLZ} et \cite{BBmoy}). Le calcul de $\overline{\Pi}(V_{k,a_p})$ a \'et\'e fait par Breuil pour $k \leq 2p+1$ (voir \cite{Br2} pour $k \leq 2p$; si $k=2p+1$ et $p\neq 2$, c'est un calcul non-publi\'e de Breuil). Voici la liste des r\'esultats que l'on obtient.

\begin{theo}\label{allred}
La r\'eduction modulo $p$ des repr\'esentations $V_{k,a_p}$ est connue dans les cas suivants.

\begin{enumerate}
\item Si $2 \leq k \leq p+1$, alors $\overline{V}_{k,a_p} = \ind(\omega_2^{k-1})$.

\item Pour $k=p+2$ :
\begin{enumerate}

\item si $1 > \val(a_p) > 0$, alors $\overline{V}_{k,a_p} = \ind(\omega_2^2)$.

\item si $\val(a_p) \geq 1$, 
et si $\lambda$ est une racine du
polyn\^ome $\lambda^2 - \overline{a_p/p} \cdot \lambda +1 = 0$,
alors : \[ \overline{V}_{k,a_p} = 
\begin{pmatrix} \omega\mu_{\lambda} & 0 \\ 
0 & \omega\mu_{\lambda^{-1}} \end{pmatrix}. \]
\end{enumerate}

\item Pour $2p \geq k \geq p+3$ :
\begin{enumerate}
\item si $1 > \val(a_p) > 0$, alors $\overline{V}_{k,a_p} =   \ind(\omega_2^{k-p})$.

\item si $\val(a_p) = 1$, et
si $\lambda=\overline{a_p/p} \cdot (k-1)$, alors :
\[ \overline{V}_{k,a_p} = 
\begin{pmatrix} \omega^{k-2} \mu_{\lambda} & 0 \\ 
0 & \omega\mu_{\lambda^{-1}} \end{pmatrix}. \]

\item si $\val(a_p) > 1$, alors $\overline{V}_{k,a_p} =   \ind(\omega_2^{k-1})$.
\end{enumerate}

\item Pour $k=2p+1$ (et $p \neq 2$) :
\begin{enumerate}
\item si $\val(a_p^2+p) < 3/2$, alors $\overline{V}_{k,a_p} =   \ind(\omega_2^2)$.
\item si $\val(a_p^2+p) \geq 3/2$, alors :
\[ \overline{V}_{k,a_p} = 
\begin{pmatrix} \omega \mu_{\lambda} & 0 \\ 
0 & \omega\mu_{\lambda^{-1}} \end{pmatrix} 
\quad\text{o\`u}\quad \lambda^2- \overline{\frac{a_p^2+p}{2 p a_p}} \cdot \lambda + 1 = 0. \]
\end{enumerate}

\item Pour $k \geq 2p+2$, les r\'esultats ne sont que partiels : 
\begin{enumerate}
\item si $\val(a_p)> \lfloor (k-2)/(p-1) \rfloor$ et si $p+1 \nmid k-1$, alors $\overline{V}_{k,a_p} =   \ind(\omega_2^{k-1})$.

\item si $\val(a_p)> \lfloor (k-2)/(p-1) \rfloor$, si $p+1 \mid k-1$ et si $i^2=-1$, alors :
\[ \overline{V}_{k,a_p} = 
\begin{pmatrix} \mu_{i} & 0 \\ 
0 & \mu_{-i} \end{pmatrix} \otimes \omega^{\frac{k-1}{p+1}}. \]
\end{enumerate}
\end{enumerate}
\end{theo}

\begin{rema}\label{forkisin}
Pour $k=p+3$ et $\val(a_p)=1$, on a un r\'esultat plus pr\'ecis (voir \cite[th\'eor\`eme 2.2.5]{BBmoy}) concernant l'extension des deux caract\`eres pour $\lambda = \pm 1$ : il existe un r\'eseau de $V_{k,a_p}$ dont la r\'eduction modulo $p$ est : 
\[\begin{pmatrix} \omega  & \ast \\ 
0 & 1 \end{pmatrix} \otimes \omega \mu_{\lambda}, \]
o\`u $\ast$ est non-trivial et peu ramifi\'e. 
\end{rema}

Pour ce qui est des repr\'esentations semi-stables, le calcul de la r\'eduction c\^ot\'e galois a \'et\'e fait dans \cite{BM2} (pour $k$ pair v\'erifiant $2 \leq k<p$) et le calcul de la r\'eduction c\^ot\'e $\G$ a \'et\'e fait dans \cite{BRst,BM} (pour $k$ pair v\'erifiant $2 \leq k \leq p+1$ et en supposant de plus que l'invariant $\mathcal{L}$ est de valuation $\geq 0$ si $k \neq 2$). Nous ne rappelons pas ici les formules, que l'on peut trouver dans \cite[\S 4.3.3]{BM2} c\^ot\'e galois et dans \cite[\S 4.5]{BRst} (pour $k=2$) et  \cite[\S 1.1]{BM} (pour $k \geq 4$) c\^ot\'e $\G$.

\vspace{20pt}
\noindent\textbf{Remerciements}: Je remercie C. Breuil pour ses encouragements tout au long de la r\'edaction de ce travail, et je le remercie de m\^eme que P. Colmez pour des discussions \'eclairantes sur plusieurs points de cet article.

\end{document}